\documentclass[12pt,twoside]{article}
\usepackage{amsmath, amsfonts, amsthm}
\usepackage{amssymb,graphics,epsfig}
\usepackage{enumerate}
\newtheorem{lemma}{Lemma}[section]
\newtheorem{theorem}{Theorem}[section]

\newtheorem{prop}{Proposition}[section]

\newcommand{\R}{ {\mathbb R} }
\newcommand{\Z}{ {\mathbb Z} }

\makeatletter \@addtoreset{equation}{section} \makeatother

\newcommand{\cqfd}{{\unskip\kern 6pt\penalty 500
\raise -2pt\hbox{\vrule\vbox to 6pt{\hrule width 6pt
\vfill\hrule}\vrule}\par}}

\newcommand{\ep}{{\varepsilon}}

\newcommand{\CCC}{C}
\newcommand{\Linf}{L^\infty}

\newcommand{\Ebar}{\overline{E}}
\newcommand{\dEbar}{\Delta\overline{E}}

\title{N-particles approximation of the Vlasov equations with
singular potential}
\author{Maxime Hauray \\ CEREMADE \\ Place du Maréchal Lattre de
Tassigny \\ 75775 Paris Cedex 16  \vspace{1cm} \\
 Pierre-Emmanuel Jabin
\\ DMA - ENS Paris \\ 45, Rue d'Ulm \\ 75005 Paris }
\date{}

\begin{document}
\maketitle

{\bf Abstract.} We prove the convergence in any time interval
of a point-particle
approximation of the
Vlasov equation by particles initially equally distributed for a
force in $1/|x|^{\alpha}$, with $\alpha \leq 1$. We introduce discrete
versions of the $L^\infty$ norm and time averages of the force field.
The core of the proof is to show
that these quantities are bounded and that consequently the minimal distance
between particles in the phase space never vanishes.


\section{Introduction}
We are interested here by the validity of the modeling of a
continuous media by a kinetic equation, with a density of presence
in phase and speed. It others words, does the many particles
follow  the evolution given by the continuous media when their
number is sufficiently large? This is a very general question and
this paper claims to give an (partial) answer only for the mean
field approach.

Let us be more precise. We study the evolution of $N$ particles,
centered at $(X_1,\dots,X_n)$ in $\R^d$ with velocities
$(V_1,\dots,V_n)$ and interacting with a central force $F(x)$. The
positions and velocities satisfy the following system of ODEs
\begin{equation} \label{ODE}
\left\{\begin{array}{l}
\dot{X_i} = V_i, \\
\dot{V_i}= E(X_i)=\displaystyle\sum_{j \neq i}
\frac{\alpha_i\,\alpha_j}{m_i}\;F(X_i-X_j),
\end{array} \right.
\end{equation}
where the initial conditions $(X_1^0,V_1^0,\dots,X_n^0,V_n^0)$ are
given. The prime example for \eqref{ODE} consists in charged particles
with charges $\alpha_i$ and masses $m_i$, in which case $F(x)=-x/|x|^3$ in
dimension three.

To easily derive from \eqref{ODE} a kinetic equation (at least
formally), it is very convenient to assume that the particles are
identical which means $\alpha_i=\alpha_j$.
Moreover we will rescale system \eqref{ODE} in time and space to work
with quantities of order one, which means
that we may assume that
\begin{equation}
\frac{\alpha_i\,\alpha_j}{m_i}= \frac{1}{N},\quad \forall i,j.
\end{equation}
We write now the Vlasov equation modelling the evolution of a
density $f$ of particles interacting with a radial force in
$F(x)$. This is a kinetic equation in the sense that the density
depends on the position and on the velocity (and of course of the
time)
\begin{equation}\begin{split}
&\partial_t f + v \cdot \nabla_x f + E(x) \cdot \nabla_v f = 0,\quad
t\in\R_+,\ x\in\R^d,\ v\in\R^d, \\
& F(x) = \nabla (\int_{x,v}
\rho(t,y)\,F(x-y)\,dy ),\\
&\rho(t,y)= \int_v f(t,x,v) \,dv.\end{split}\label{vlasov}
\end{equation}
Here  $\rho$ is the spatial density and the initial density $f^0$
is given.

When the number $N$ of particles is large, it is obviously easier
to study (or solve numerically) \eqref{vlasov} than \eqref{ODE}.
Therefore it is a crucial point to determine whether
\eqref{vlasov} can be seen as a limit of \eqref{ODE}.

Remark that if $(X_1,\ldots,X_N,V_1,\ldots,V_n)$ is a solution of
\eqref{ODE}, then the measure
$$
\mu_N(t)= \frac{1}{N} \sum_{i=1}^n \delta(x-X_i(t))\otimes\delta
(v-V_i(t))
$$
is a solution of the Vlasov equation in the sense of
distributions. And the question is whether a weak limit $f$ of
$\mu_N$ solves \eqref{vlasov} or not. If $F$ is smooth, then it is
indeed the case as it is proved in the book by Spohn \cite{Sp}.
The purpose of this paper is to justify this limit if
\begin{equation}
|F(x)|\leq \frac{C}{|x|^\alpha},\quad |\nabla F(x)|\leq
\frac{C}{|x|^{1+\alpha}}\quad |\nabla^2 F(x)|\leq
\frac{C}{|x|^{2+\alpha}},\quad \forall x\neq 0,
\label{boundF}\end{equation}
for $\alpha<1$, which is the first rigorous proof of the limit in a
case where $F$ is not necessarily bounded.

Before being more precise concerning our result, let us explain
what is the meaning of \eqref{ODE} in view of the singularity in
$F$. Here we assume either that we restrict ourselves to the
initial configurations for which there are no collisions between
particles over a time interval $[0,\ T]$ with a fixed $T$,
independent of $N$. Or we assume that $F$ is regular or
regularized but that the norm $\|F\|_{W^{1,\infty}}$ may depend on
$N$; This procedure is well presented in \cite{Bat} and it is the
usual one in numerical simulations (see \cite{VA} and \cite{Wo}).
In both cases, we have classical solutions to \eqref{ODE} but the
only bound we may use is \eqref{boundF}.

Other possible approaches would consist in justifying that the set of
initial configurations $X_1(0),\ldots,X_N(0), V_1(0),\ldots, V_N(0)$
for which there is at least one collision,
is negligible or that it is possible to define a solution (unique or
not) to the dynamics even with collisions.

Finally notice that the condition $\alpha<1$ is not unphysical.
Indeed if $F$ derives from a potential, $\alpha=1$ is the critical
exponent for which repulsive and attractive forces seem very
different. In other words, this is the point where the behavior of
the force when two particles are very close takes all its
importance.
%
%
\subsection{Important quantities \label{defquant}}
%
The derivation of the limit requires a control on many
quantities. Although some of them are important only at the discrete
level, many were already used to get the existence of strong solutions
to the Vlasov-Poisson equation (we refer to \cite{Ho1}, \cite{Ho2} and
\cite{Pf}, \cite{Sc} as being the closest from our method).

The first two are quite natural and are bounds on the size of the
support of the initial data in space and velocity, namely we introduce
\begin{equation}
R(T)=\sup_{t \in [0,T], \,i=1,\dots N} X_i( t),\quad
 K(T)=\sup_{t \in [0,T], \, i=1,\dots N} V_i(
t).\label{defRK}
\end{equation}

Of course $R$ is trivially controlled by $K$ since
\begin{equation}
R(T)\leq R(0)+T\,K(T).
\label{RV}
\end{equation}
Now a very important and new parameter is the discrete scale of
the problem denoted $\ep$. This quantity represents roughly the
minimal distance between two particles or the minimal time
interval which the discrete dynamics can see. We fix this
parameter from the beginning and somehow the main part of our work
is to show that it is indeed correct, so take
\begin{equation}
\ep=\frac{R(0)}{N^{1/2d}}.
\label{defep}
\end{equation}
At the initial time, we will choose our approximation  so that the
minimal distance between two particles will be of order $\ep$.

The force term cannot be bounded at every time for the discrete
dynamics (a quantity like $F \star \rho_N$ is not bounded even in
the case of  free transport), but we can expect that its average
on a short interval of time will be bounded. So we denote
\begin{equation}
\Ebar(T) = \sup_{t \in [t_0,T-\ep],i=1,\dots,N} \left\{
\frac{1}{\ep} \int_t^{t+\ep} |E(X_i(s))| \,ds \right\},
\label{defEbar}
\end{equation}
with for $T< \ep$
\begin{equation}
\Ebar(T) = \sup_{i=1,\dots,N} \left\{ \frac{1}{\ep}
\int_0^{T} |E(X_i(s))| \,ds \right\},
\label{defEbar'}
\end{equation}
thus obtaining a continuous definition. Moreover we denote by
$E^0$ the supremum over all $i$ of $|E(X_i(0))|$.

This definition comes from the following intuition. The force is
big when two particles are close together. But if their speeds are
different, they won't stay close a long time. So we can expect the
interaction force between these two particles to be integrable in
time even if they "collide". They just remain the case of two
close particles with almost the same speed. To estimate the force
created by them, we need an estimate on their number. On way of
obtaining it is to have a bound on
\begin{equation}
m(T) = \sup_{t \in [0,T], i \neq j}
\frac{\ep}{|X_i(t)-X_j(t)|+|V_i(t)-V_j(t)|}, \label{defm}
\end{equation}
The control on $m$ requires the use of a discretized derivative of
$E$, more precisely we define for any exponent
$\beta\in\left.\right]1,\ d-\alpha\left[\right.$ which we note
also satisfies $\beta<2d-3\alpha$ ($\beta=1$ would be
enough for short time estimates)
\begin{equation}
 \dEbar(T) = \sup_{t \in
[t_0,T-\ep]}\,\sup_{i,j=1,\dots,N,}
\left\{ \frac{1}{\ep}
\int_t^{t+\ep} \frac{|E(X_i(s))-E(X_j(s))|}{\ep^\beta+|X_i(s)-X_j(s)|}
\,ds \right\},
\label{defdEbar}
\end{equation}
with as for $\Ebar$, when $T<\ep$
\begin{equation}
 \dEbar(T) = \sup_{i,j=1,\dots,N}
\left\{ \frac{1}{\ep}
\int_0^{T} \frac{|E(X_i(s))-E(X_j(s))|}{\ep^\beta
+|X_i(s)-X_j(s)|} \,ds \right\}.
\label{defdEbar'}
\end{equation}
 Now, we introduce what we called the discrete infinite norm of
 the distribution of the particle $\mu_N$. This quantities is the
 supremum over all the boxes of size $\ep$ of the total mass they
 contains divided by the size of the box. That is, for a measure $\mu$
we denote
\begin{equation}
 \|\mu\|_{\infty,\ep}= \frac{1}{(2\ep)^6} \sup_{(x,v) \in \R^6}
\left\{ \mu(B_{\infty}((x,v),\ep)) \right\}.
\label{deflinfty}
\end{equation}
where
$B_{\infty}((x,v),\ep)$ is the ball of radius $\ep$ centered at
$(x,v)$ for the infinite norm. Note that we may bound
$\|\mu_N(T,\cdot)\|_{\infty,\ep}$ by
\begin{equation}
\|\mu_N(T,\cdot)\|_{\infty,\ep} \leq \left(4\,m(T)\right)^{2d}.
\label{mlinf}
\end{equation}

All the previous quantities (except for $\ep$) will always be
assumed to be bounded at the initial time $T=0$ uniformly in $N$.
\subsection{Main results}
The main point in the derivation of the Vlasov equation is to obtain a
control on the previous quantities. We first do it for a short time as
given by
\begin{theorem} If $\alpha<1$,
there exists a time $T$ and a constant $c$ depending only on $R(0)$, $K(0)$,
$m(0)$  but not on $N$ such that for some $\alpha<\alpha'<3$
\[ \begin{split}
&R(T)\leq 2\; (1+R(0)),\quad K(T)\leq 2\; (1+K(0)),\quad m(T)\leq 2\,m(0),\\
&\Ebar(T)\leq
c\;(m(0))^{2\alpha'}\,(K(0))^{\alpha'}\,
(R(0))^{\alpha'-\alpha}
,\quad \sup_{t\leq T}
\|\mu_N(t,\cdot)\|_{\infty,\ep}\leq (8\,m(0))^{2d}.\\
\end{split}
\]
\label{bornes}
\end{theorem}
{\bf Remark}\par The constant $2$ is of course only a matter of
convenience. But the choice of a different constant is not really
helpful. This result is valid only for a short time in essence and
increasing the chosen constant for instance, increases only
slightly the time $T$ who in any case cannot pass a critical
value.

\medskip

This theorem can be extended on any time interval
\begin{theorem} For any time $T>0$, there exists a function
$\tilde N$ of $R(0)$, $K(0)$, $m(0)$, and $T$
and a constant $C(R(0),K(0),m(0),T)$ such
that if $N\geq \tilde N$ then
\[
R(T),\ K(T),\ m(T),\ \Ebar(T)\leq C(R(0),K(0),m(0),T).
\]
\label{longtimeexist}
\end{theorem}
>From this last theorem, it is easy to deduce the main result of this
paper which reads
\begin{theorem} Consider a time $T$ and
sequence $\mu_N(t)$ corresponding to solutions
to \eqref{ODE} such that $R(0)$, $K(0)$ and $m(0)$ are bounded
uniformly in $N$. Then any weak limit $f$ of $\mu_N(t)$ in
$L^\infty([0,\ T],\ M^1(\R^{2d}))$ belongs to $L^\infty([0,\ T],\
L^1\cap L^\infty(\R^{2d}))$, has compact support and is a solution to
\eqref{vlasov}. \label{derivation}
\end{theorem}
Of course the main limitation of our results is the condition
$\alpha<1$ and the main open question is to know what happens when
$\alpha\geq 1$. However this condition is not only technical and new
ideas will be needed to prove something for $\alpha\geq 1$.

It would also be interesting to extend our result to more
complicated forces like the ones found in the formal derivation of
\cite{JP}.

The derivation of kinetic equations is an important question both
for numerical and theoretical aspects. We already mentioned the
works of Batt \cite{Bat}, Spohn \cite{Sp}, Victory and Allen
\cite{VA} and Wollmann \cite{Wo} for Vlasov equations. Another
interesting case concerns Boltzmann equation, for which we refer
to the book by Cercignani, Illner and Pulvirenti \cite{CIP} and
the paper by Illner and Pulvirenti \cite{IP}.

On the other hand, the derivation of macroscopic equations is
usually easier and some results are already known (although not
since a very long time) even in cases with singularity. In
particular and that is more or less the macroscopic equivalent of
our result, the convergence of the point vortex method for $2-D$
Euler equations was obtained by Goodman, Hou and Lowengrub
\cite{GHL} (see also the works by Schochet \cite{Sch1} and
\cite{Sch2}). The main difficulty for macroscopic systems is to
control the minimal distance between two particles (which is not
possible in the kinetic framework) as it is also clear in
\cite{JO}.

Our method of proof makes full use of the characteristics and of
the procedures developed to get for the Vlasov-Poisson equation in
dimension two and three. This method was introduced by Horst in
\cite{Ho1} and \cite{Ho2} and was successfully used to prove the
existence of strong solutions in large time in \cite{Pf} and
\cite{Sc} and at the same time by Lions and Perthame in \cite{LP2}
using the moments (see also \cite{GJP} for a slightly simpler
proof and \cite{Pe2} for an application to the asymptotic behavior
of the equation). These results were extended to the periodic case
by Batt and Rein in \cite{BR} and to the
Vlasov-Poisson-Fokker-Planck equation by Bouchut in \cite{Bo}. In
particular the necessity to integrate in time to control the
oscillations of the force also appears in the proof of $L^\infty$
bounds for the Vlasov-Poisson-Fokker-Planck equation by Pulvirenti
and Simeoni in \cite{PuSi}. We refer to the book by Glassey
\cite{Gl} for a general discussion of the existence theory for
kinetic equations.

In the rest of the paper, $C$ will denote a generic constant,
depending maybe on $R(0)$, $K(0)$, or
$m(0)$  but not on $N$ or any other quantity. We first prove Theorem
\ref{bornes}, then we show a preservation of discrete $L^\infty$ norms
which proves Theorem \ref{longtimeexist}. In the last section we
explain how to deduce Theorem \ref{derivation}.
%
\section{Proof of Theorem \ref{bornes}}
The first steps are to estimate all quantities in term of
themselves. Then if this is done correctly it is possible to deduce
bounds for them on a short interval of time.

\subsection{Estimate on $\Ebar$}
We are in fact able to estimate any time average of the force field,
more precisely
\begin{lemma} For any $\alpha'$ with $\alpha<\alpha'<3$, assume that
\[
m(t_0)\leq \frac{1}{12\,\ep\,K(t_0)\,\dEbar(t_0)},
\]
then the following
inequality holds
\[ \begin{split}
\Ebar(t_0) \leq C\, ( &\|\mu_N\|_{\infty,\ep}^{\alpha'/d}\,
K^{\alpha'}\,R^{\alpha'-\alpha} + \ep^{d-\alpha}\, \|\mu_N\|_{\infty,\ep}\,
K^{2d-\alpha} \\
&\ + \ep^{2d-3\alpha}\,
\|\mu_N\|_{\infty,\ep}\, K^{d-\alpha} \Ebar^d\,K^d),
\end{split}
\]
 where we use the values of $\|\mu_N\|_{\infty,\ep}$, $R$, $K$, $m$
and $\Ebar$ at
the time $t_0$.
\label{moytps}
\end{lemma}
Of course if any of the above quantity is infinite then the result is
obvious. This lemma could appear stupid since we control
$\Ebar(t_0)$ by itself (and with a power larger than $1$ in
addition). But the point is that except for the first term, the other
two are very small because of the $\ep$ in front of them so that they
almost do not count.

\begin{proof} For any index $i$, we bound for any $t_1$ less than
$t_0-\ep$ or $t_1=0$ if $t_0<\ep$
\[
I^i=\frac{1}{\ep}\int_{t_1}^{t_1+\ep} |E(X_i(s))|\,ds,
\]
which will give the desired result by taking the supremum over all
$i$. So now let us fix $t_1$ and $i$ (we choose $i=1$ for
simplicity).

 We introduce the following decomposition among the particles:
We define
\begin{equation}
C_k=\biggl\{i\,\Bigl|\;
3\,\ep\,K(t_0)\,
2^{k-1}<|X_i(t_1)-X_1(t_1)|\leq 3\,\ep\,K(t_0)\, 2^{k}\biggr\}.
\end{equation}
Therefore the index
$k$ varies from $1$ to $k_0=(\ln (R/4\,\ep\,K(t_0)))/\ln 2$. And we
denote by
$C_0$ the rest, that is the set of indices $i$ such that
$|X_i(t_1)-X_1(t_1)|\leq 3\,\ep\,K(t_0)$.

Consequently we decompose the term $I^1$  into two parts,
$I_1 + I_{C_0}$ with
\begin{equation}
 I_1 = \sum_{k=2}^{k_0} \sum_{i \in C_k}
\frac{1}{\ep}\int_{t_1}^{t_0}
\frac{1}{N\,|X_1(t)-X_i(t)|^{\alpha}} \,dt,
\label{defI1}
\end{equation}
and
$I_{C_0}$ is the sum of the same terms over the particles of
$C_0$.

\medskip

{\em Step 1: Stability of the $C_k$.} Given their definition, the
$C_k$ enjoy the following property, for any $i\in C_k$ with $k>1$, we
have for any $t\in [t_1,\ t_0]$
\[
|X_1(t)-X_i(t)|\geq \ep\,K(t_0)\,2^{k-1}.
\]
Indeed, we of course know that
\[
\left|\frac{d}{dt}(X_i(t)-X_1(t))\right|=|V_i(t)-V_1(t)|\leq 2\,K(t_0),
\]
and then
\[\begin{split}
|X_1(t)-X_i(t)|&\geq |X_1(t_1)-X_i(t_1)|-2\,(t_0-t_1)\,K(t_0)\\
& \geq
3\,\ep\,K(t_0)\, 2^{k-1}-2\,\ep\,K(t_0),
\end{split}\]
with the corresponding result since $k\geq 1$. Of course the same
argument also shows that if $i\in C_0$ then for any $t\in[t_1,\ t_0]$,
\[
|X_1(t)-X_i(t)|\leq 5\,\ep\,K(t_0).
\]

\medskip

{\em Step 2: Control of $I_1$.} Using the result from the previous
step, we deduce that for any $i\in C_k$ with $K\geq 1$,
\[
\frac{1}{|X_i(t)-X_1(t)|^\alpha}\leq \frac{C\,2^{-\alpha
k}}{\ep^\alpha\,(K(t_0))^\alpha}.
\]
On the other hand, we have of course $|C_k|\leq N$ and
moreover $|C_k|\leq C\, \ep^{-2d}\, K^{2d}\, \ep^d\, 2^{d\,
k}\times \|\mu_N\|_{\infty,\ep}$ according to the very definition of
this discrete $L^\infty$ norm \eqref{deflinfty}.
Consequently for any $\alpha'<d$, since $\ep^{2d}=C/N$,
interpolating between these two values, we get
\[
|C_k|\leq C\,N\, (K(t_0))^{2\alpha'}\,
\ep^{\alpha'}\, 2^{\alpha'k}\times
\|\mu_N(t_0,.)\|_{\infty,\ep}^{\alpha'/d}.
\]
Using these last two bounds in \eqref{defI1} and summing up, we obtain
\[\begin{split}
I_1&\leq \sum_{k=1}^{k_0} |C_k|\times N^{-1}\;
(K(t_0))^{-\alpha}\,\ep^{-\alpha}\, 2^{-\alpha\,k}\\
&\leq C\,
\|\mu_N\|_{\infty,\ep}^{\alpha'/d}\, K^{2\alpha'-\alpha}\,
\ep^{\alpha'-\alpha}
\sum_{k=1}^{k_0} 2^{(\alpha'-\alpha)k}.
\end{split}\]
Eventually for any $\alpha<\alpha'<d$, we deduce that
\begin{equation}\begin{split}
I_1\leq C\,  \|\mu_N\|_{\infty,\ep}^{\alpha'/d}\,
K^{2\alpha'-\alpha}\,\ep^{\alpha'-\alpha}\;
2^{(\alpha'-\alpha)k_0}\leq C\, \|\mu_N\|_{\infty,\ep}^{\alpha'/d}\,
R^{\alpha'-\alpha}\, K^{\alpha'},
\end{split}
\label{estE1}
\end{equation}
the values being taken at $t_0$,
which gives the first term in Lemma \ref{moytps}. Before dealing with
the remaining term, we point out that here we have never used the
condition $\alpha<1$ and that for this term
the same computation would be valid for
any $\alpha<3$.

\medskip

{\em Step 3: Redecomposition of $C_0$.}
  For the force induced by the particles in $C_0$, we divided again
this set into several parts.
Set
\begin{equation}
Q_l = \biggl\{j \in C_0\, \Bigl|\; 3\,\ep\,\Ebar(t_0)\,2^{l-1} \leq
  |V_1(t_1)-V_j(t_1)| \leq 3\,\ep\, \Ebar(t_0)\,2^l\biggr\} ,
\label{defQl}
\end{equation}
for $l \geq 1$.
Remark that $Q_l = \emptyset$ if $l > l_0 = \ln
(K(t_0)/(\ep\, \Ebar(t_0)))/\ln 2$. Therefore the rest $Q_0$ is
defined by
\[
Q_0 = \biggl\{j \in C_0 \Bigl| |V_1(t_1)-V_j(t_1)| \leq
3\,\ep\,\Ebar(t_0)\biggr\}.
\]
As before we decompose $I_{C_0}$ in a sum of $I_2$ and a remainder
$I_{Q_0}$ with
\begin{equation}
I_2=\sum_{l=1}^{l_0}\sum_{j\in
Q_l}\frac{1}{\ep}\int_{t_1}^{t_0}\frac{dt}{N\,|X_j(t)-X_1(t)|},
\label{defI2}
\end{equation}
and for $I_{Q_0}$ the same sum but on the indices $j\in Q_0$ of course.

The idea behind this new decomposition is that although the
particles in $Q_l$ witth $l\geq 1$ are close to $X_1$,
their speed is different from $V_1$. So even if
  they come very close to $X_1$ they will stay close only for a very
short time.
  Since the singularity of the potential is not too high, we
  will be able to bound the force.

\medskip

{\em Step 4: Stability of the $Q_l$.} Just as for the $C_k$, we may
prove that for any time $t$ in $[t_1,\ t_0]$ and any $j\in A_l$ with
$l\geq 1$
\[
|V_j(t)-V_1(t)|> \ep\,\Ebar(t_0)\,2^{l-1}.
\]
This is again due to the fact that
\[
|V_j(t)-V_j(t_1)|\leq \int_{t_1}^{t_0} |E(X_j(s))|\,ds\leq
\ep\,\Ebar(t_0),
\]
so that in fact the result is even more precise in the sense that the
relative velocity $V_j(t)-V_1(t)$ remains close to $V_j(t_1)-V_1(t_1)$
up to exactly $\ep\,\Ebar(t_0)$.

\medskip

{\em Step 5: Control of $I_2$.} Given this previous point, for any
$j\in Q_l$ with $l>0$ and any $t\in[t_1,\ t_2]$, we have, denoting by
$t_m$ the time in the interval $[t_1,\ t_0]$ where $|X_j(t)-X_1(t)|$
is minimal
\[
|X_1(t) - X_j(t)| \geq \left| |X_1(t_m) - X_j(t _m)| - \frac{1}{2}
(t-t_m)|V_1(t_m)-V_j(t_m)| \right|.
\]

Then,
\begin{eqnarray*} \frac{1}{\ep} \int_{t_1}^{t_0}
\frac{1}{|X_1(t) - X_j(t)|^\alpha}  \,dt & \leq  & \frac{C}{\ep}\;
|V_1(t_m)-V_j(t_m)|^{-\alpha} \ep^{1-\alpha}\\
  & \leq & C \,\ep^{-2\alpha}\,(\Ebar(t_0))^{-\alpha}\,2^{-\alpha l}.
\end{eqnarray*}
Summing up on $l$, we obtain
$$ |I_2| \leq C\sum_{l=1}^{l_0} |Q_l|
\frac{1}{N}\,\ep^{-2\alpha}\, (\Ebar(t_0))^{-\alpha}\,
2^{-\alpha l}.
$$
We bound $|Q_l|$ by $|Q_l| \leq C\,\|\mu\|_{\infty,\ep}\,(K(t_0)\,\ep)^d
  (2^l\,\Ebar(t_0)\,\ep)^d$ using again the definition of the discrete
$L^\infty$ norm. It gives us
  \begin{eqnarray*}
  |I_2| & \leq & C\, (K(t_0))^d\, (\Ebar(t_0))^{d-\alpha}\, \ep^{2d-2\alpha}\,
\|\mu_N\|_{\infty,\ep}\times
  \sum_{l=2}^{l_0} 2^{(d-\alpha)l} \\
   & \leq & C\, (K(t_0))^d\, (\Ebar(t_0))^{d-\alpha}\,
\|\mu_N\|_{\infty,\ep}\,
 \ep^{2d-2\alpha}\,\left(\frac{K(t_0)}{\Ebar(t_0)\,
   \ep}\right)^{d-\alpha}\\
   & \leq & C\, (K(t_0))^{2d-\alpha}\, \|\mu_N\|_{\infty,\ep}\,
\ep^{d-\alpha},
  \end{eqnarray*}
which is indeed the second term in Lemma \ref{moytps}.

\medskip

{\em Step 6: Control on $I_{Q_0}$.} The first point to note is that
for any $j\in Q_0$ and any $t\in[t_1,\ t_0]$
by the definition \eqref{defdEbar} of $\dEbar$ and the stability of $C_0$
\[
|V_j(t)-V_1(t)-V_j(t_1)-V_1(t_1)|\leq 5\, \ep^2\, K(t_0)\,\dEbar(t_0).
\]
Consequently it is logical to decompose (again) $Q_0$ in $Q_0'\cup
Q_0''$ and $I_{Q_0}$ in $I_{Q_0'}+I_{Q_0''}$ with
\[
Q_0'=\biggl\{j\in Q_0\,\Bigl|\ |V_j(t_1)-V_1(t_1)|\geq 6\,\ep^2\,
K(t_0)\, \dEbar(t_0)\biggr\},
\]
and $I_{Q_0'}$, $I_{Q_0''}$ the sums on the corresponding indices.

Then for any $j\in Q_0'$, the same computation as in the fifth step,
shows that
\[
\frac{1}{\ep}\int_{t_1}^{t_0}\frac{dt}{N\,|X_j(t)-X_1(t)|^\alpha}\leq
C\,\ep^{2d-3\alpha}\,(K(t_0))^{-\alpha}\,(\dEbar(t_0))^{-\alpha}.
\]
But of course the cardinal of $Q_0'$ is bounded by the one of $Q_0$
and using as always the discrete $L^\infty$ bound
\[
|Q_0'|\leq C\, (K(t_0))^d\,(\Ebar(t_0))^d\,\|\mu_N\|_{\infty,\ep}.
\]
Eventually that gives
\[
I_{Q_0'}\leq
C\,\ep^{2d-3\alpha}\,(K(t_0))^{d-\alpha}\,(\Ebar(t_0))^d\,
\|\mu_N(t_0,.)\|_{\infty,\ep},
\]
since $\dEbar(t_0)$ comes with a negative exponent and being non
decreasing, may be bounded by the initial value. So $I_{Q_0'}$
corresponds to the third term in the lemma.

Let us conclude the proof with the bound on $I_{Q_0''}$. Of course if
$j\in Q_0''$ then for any $t\in[t_1,\ t_0]$,
\[
|V_j(t)-V_1(t)|\leq 11\,\ep^2\,K(t_0)\,\dEbar(t_0).
\]
Now we use the definition \eqref{defm} of $m$ and the assumption in
the lemma to deduce that
\[
|X_j(t)-X_1(t)|\geq \frac{\ep}{m(t_0)}-|V_j(t)-V_1(t)|\geq \ep^2\,K(t_0)\,
\dEbar(t_0).
\]
We bound $|Q_0''|$ by $Q_0$ which is the best we can do since the
discrete $L^\infty$ norm cannot see the scales smaller than $\ep$ and
we obtain
\[
I_{Q_0''}\leq C\, \ep^{2d-2\alpha}\,
(K(t_0))^{d-\alpha}\,(\Ebar(t_0))^d\,\|\mu_N(t_0,.)\|_{\infty,\ep},
\]
which is dominated by the previous term and the third term in the lemma.

Before ending the proof
we wish to note that the condition $\alpha<1$ was only used to get the
second term in the fifth step and in the last step
and the bound on $m$ was only required
in this last step.
\end{proof}
\subsection{Estimate on $\dEbar$}
We may show the following with the same remarks as for Lemma \ref{moytps},
 \begin{lemma} For any $\alpha'$ with $\alpha<\alpha'<3$, assume that
\[
m(t_0)\leq \frac{1}{12\,\ep\,K(t_0)\,\dEbar(t_0)},
\]
then the following
inequality holds
\[ \begin{split}
\dEbar(t_0) \leq C\, ( &\|\mu_N\|_{\infty,\ep}^{(1+\alpha')/d}\,
K^{1+\alpha'}\,R^{\alpha'-\alpha} + \ep^{d-\alpha-\beta}\,
\|\mu_N\|_{\infty,\ep}\,
K^{2d-\alpha} \\
&\ + \ep^{2d-3\alpha-\beta}\,
\|\mu_N\|_{\infty,\ep}\, K^{d-\alpha} \Ebar^d\,K^d),
\end{split}
\]
 where we use the values of $\|\mu_N\|_{\infty,\ep}$, $R$, $K$, $m$
and $\Ebar$ at
the time $t_0$.
\label{moytpsdE}
\end{lemma}
\begin{proof} The proof follows the same procedure as for Lemma
\ref{moytps} with exactly the same decompositions. We have to bound,
since as before the choice of the indices does not matter
\[
\Delta I=\frac{1}{\ep} \int_{t_1}^{t_0}
\frac{|E(X_1(t))-E(X_2(t))|}{\ep^\beta+|X_1(t)-X_i(t)|} \,dt.
\]
We introduce the same decomposition as for the proof of Lemma
\ref{moytps} except that now we have two decompositions: One around
$X_1$ denoted by a superscript $1$ and another one around $X_2$. So
we set for $\gamma=1,\;2$
\[
C_k^\gamma=\biggl\{i\,\Bigl|\;
3\,\ep\,K(t_0)\,
2^{k-1}<|X_i(t_1)-X_\gamma(t_1)|\leq 3\,\ep\,K(t_0)\, 2^{k}\biggr\},
\]
with $C_0^1$ and $C_0^2$ the corresponding remaining indices. We
also denote $C_0=C_0^1\cup C_0^2$.

Then of course
\[\begin{split}
\Delta I&=\sum_{i\neq C_0} \frac{1}{\ep}
\!\int_{t_1}^{t_0}\!\!\left|F(X_1-X_i)-F(X_2-X_i)\right|\times
\frac{dt}{\ep^\beta+|X_1(t)-X_i(t)|}+\sum_{i\in C_0}\ldots\\
 &=\Delta I_1+\Delta I_{C_0}.\\
\end{split}\]
The second term is easier to bound as we simply write
\[
\Delta I_2\leq \sum_{i\in C_0}\frac{1}{\ep^{1+\beta}}
\int_{t_1}^{t_0}\! \left(\frac{1}{N\,|X_i(t)-X_1(t)|^\alpha}
+\frac{1}{N\,|X_i(t)-X_2(t)|^\alpha}\right)\;dt.
\]
Now we do exactly what we did for $I_{C_0}$ in the proof of Lemma
\ref{moytps} and we get as a bound the two last terms in the estimate
for $\Ebar$ in this lemma, divided by $\ep^\beta$, which are exactly the two
last term in the estimate for $\dEbar$ in Lemma \ref{moytpsdE}. Note
by the way
that the terms where we sum on $C_0^1$ for instance but with $X_i-X_2$
in the denominator are in fact even easier to handle.

\medskip

Now for $\Delta I_1$, we observe that for $i\neq C_0$ then for any $t$
\[\begin{split}
|F(X_1(t)-X_i(t))-&F(X_2(t)-X_i(t))|\leq C |X_1(t)-X_2(t)|\\
& \times
\left(\frac{1}{N\,|X_1(t)-X_i(t)|^{\alpha+1}}
+\frac{1}{N\,|X_1(t)-X_i(t)|^{\alpha+1}}\right),
\end{split}\]
since it is always possible to find a regular path $x_t(s)$ of length
less than $2\,|X_1(t)-X_2(t)|$ such that
$x_t(0)=X_1(t)$, $x_t(1)=X_2(t)$ and $|x_t(s)-X_i(t)|$ is always
larger than the minimum between $|X_1(t)-X_i(t)|$ and
$|X_2(t)-X_i(t)|$.

Therefore we are led to make exactly the same computation as for $I_1$
with $\alpha+1$ instead of $\alpha$ which gives the desired
result. This is possible only because in the estimate on $I_1$, we
never use the condition $\alpha<1$.
\end{proof}

\subsection{Control on $m$ and $K$}
We prove the
\begin{lemma} Assume that
\[
m(t)\leq \frac{1}{\ep^{\beta-1}},
\]
then we also have that
\[
m(t)\leq m(0)\times e^{Ct+C\ep\,\dEbar(t)+C\int_0^t \dEbar(s)\,ds},
\]
and we may eliminate the $\ep\,\dEbar(t)$ term if $t>\ep$.
\label{boundm}
\end{lemma}
Note that we still need an assumption on $m$ but it is a bit different
(and somewhat ``harder'' to satisfy) than the corresponding one for
Lemmas \ref{moytps} and \ref{moytpsdE}.
\begin{proof} We consider any two indices $i\neq j$. Then we write
\[\begin{split}
\frac{d}{ds}&\left(\frac{\ep}{|X_i(s)\!-\!X_j(s)|+|V_i(s)\!-\!V_j(s)|}\right)=
\frac{\ep}{(|X_i(s)\!-\!X_j(s)|+|V_i(s)\!-\!V_j(s)|)^2}\\
&\times \Bigl( \frac{X_i-X_j}{|X_i-X_j|}\cdot (V_i-V_j)
+\frac{V_i-V_j}{|V_i-V_j|}\cdot (E(X_i)-E(X_j))\Bigr)\\
&\leq\frac{\ep\,(|V_i(s)-V_j(s)|
+|E(X_i(s))-E(X_j(s))|)}{(|X_i(s)\!-\!X_j(s)|
+|V_i(s)\!-\!V_j(s)|)^2}.
\end{split}
\]
Since $m(t)\leq \ep^{1-\beta}$, the same is true of $m(s)$ and at
least one of
the quantities $|X_i(s)-X_j(s)|$ and $|V_i(s)-V_j(s)|$ is larger than
$\ep^\beta/2$, therefore
\[\begin{split}
\frac{d}{ds}\left(\frac{\ep}{|X_i(s)\!-\!X_j(s)|+|V_i(s)\!-\!V_j(s)|}\right)
&\leq \frac{C\,\ep}{|X_i(s)\!-\!X_j(s)|+|V_i(s)\!-\!V_j(s)|}\\
&\times \left(1+\frac{|E(X_i(s))-E_(X_j(s))|}{\ep^\beta
+|X_i(s)-X_j(s)|}\right).
\end{split}
\]
But by the definition of $\dEbar$, see \eqref{defdEbar}, we know that
for $t>\ep$
\[
\int_\ep^t \frac{|E(X_i(s))-E_(X_j(s))|}{\ep^\beta
+|X_i(s)-X_j(s)|}\;ds\leq \,\int_0^t \dEbar(s)\,ds,
\]
and of course for $t<\ep$
\[
\int_0^t  \frac{|E(X_i(s))-E_(X_j(s))|}{\ep^\beta
+|X_i(s)-X_j(s)|}\;ds\leq \ep \dEbar(t).
\]
Hence, integrating in time, we find
\[\begin{split}
\frac{\ep}{|X_i(s)\!-\!X_j(s)|+|V_i(s)\!-\!V_j(s)|}&\leq
\frac{\ep}{|X_i(0)\!-\!X_j(0)|+|V_i(0)\!-\!V_j(0)|}\\
&\times e^{Ct+C\ep\,\dEbar(t)+C\int_0^t \dEbar(s)\,ds},
\end{split}\]
wich after taking the supremum in $i$ and $j$ is precisely the lemma.
\end{proof}

\medskip

As for $K$, using the equation that $\dot V_i(t)=E_i(X_i(t))$, we may
prove by the same method which we do not repeat, the result
\begin{lemma}
We have that for any $t$
\[
K(t)\leq K(0)+Ct+C\ep\,\Ebar(t)+C\int_0^t \Ebar(s)\,ds.
\label{boundK}
\]
\end{lemma}
\subsection{Conclusion on the proof of Theorem \ref{bornes}}
Here (but only in this subsection)
for a question of clarity, we keep the notation $C$ for the
constants appearing in Lemmas \ref{moytps}, \ref{moytpsdE},
\ref{boundm} and \ref{boundK} and we denote by $\tilde C$ any other
constant depending only on $R(0)$, $K(0)$ and $m(0)$.

We assume that on a time interval $[0,\
T]$, we have (for a given $\alpha'$)
\begin{equation}\begin{split}
&m(t)\leq 2\,m(0),\quad \Ebar(t)\leq
2C\,2^{8\alpha'-\alpha}\,(m(0))^{2\alpha'}\,
(K(0))^{\alpha'}\,(R(0))^{\alpha'-\alpha},\\
&K(t)\leq 2\,(1+K(0)),\quad R(0)\leq
2\,(1+R(0)),\quad \forall\,t\in\,[0,\ T],
\end{split}
\label{shorttimeas}
\end{equation}
which we may always do since all these quantities are continuous in
time (although they may a priori increase very fast).

Then we show that if $T$ is too small we have in fact the same
inequalities but with a $3/2$ constant instead of $2$. By
contradiction this of course shows that we can bound $T$ from below in
terms of only $R(0)$, $K(0)$ and $m(0)$ and it proves Theorem
\ref{bornes} with $c=C\times 2^{8\alpha'-\alpha+1}$.

\medskip

First of all, we note that since $m(t)\leq 2 m(0)$, we may apply Lemmas
\ref{moytps}, \ref{moytpsdE}, and \ref{boundm}. Furthermore we
immediately know from \eqref{mlinf} that
\[
\|\mu_N(t,.)\|_{\infty,\ep}\leq  (8\,m(0))^{2d}.
\]

Let us start with
Lemma \ref{moytps}, using the assumption \eqref{shorttimeas} we
deduce that for any $t\in[0,\ T]$,
\[
\Ebar(t)\leq C\, 2^{8\alpha'-\alpha}\,
(m(0))^{2\alpha'}\,(K(0))^{\alpha'}\, (R(0))^{\alpha'-\alpha}
+\tilde C\,\ep^{d-a}+\tilde C\,\ep^{2d-3\alpha}.
\]
For $\ep$ small enough this proves that
\[
\Ebar(t)\leq \frac{3\,C}{2}\, 2^{8\alpha'-\alpha}\,
(m(0))^{2\alpha'}\,(K(0))^{\alpha'}\, (R(0))^{\alpha'-\alpha},
\]
which is the first point.

\medskip

Next applying Lemma \ref{moytpsdE}, we deduce that for any $t\in[0,\ T]$
\[
\dEbar(t)\leq \tilde C.
\]
From Lemma \ref{boundm}, we obtain that
\[
m(t)\leq m(0)\times e^{\tilde C T},
\]
so if $T$ is such that $\tilde C\,T<\ln(3/2)$ then we get
\[
m(t)\leq \frac{3}{2}\,m(0).
\]

\medskip

Lemma \ref{boundK} implies that for $t\in[0,\ T]$
\[
K(t)\leq K(0)+\tilde C\,T,
\]
so that again for $T$ small enough
\[
K(t)\leq \frac{3}{2}\,(1+K(0)).
\]

\medskip

Eventually thanks to relation \eqref{RV}, we know that for $t\in[0,\
T]$
\[
R(t)\leq R(0)+T\,K(t)\leq R(0)+\tilde C\,T,
\]
hence the corresponding estimate for $R$ provided $\tilde C\,T\leq
3/2$.

In conclusion we have shown that if \eqref{shorttimeas} holds and if
$T$ is smaller than a given time depending only on $R(0)$, $K(0)$ and
$m(0)$ then the same inequalities are true with $3/2$ instead of
$2$. By the continuity of $R$, $K$, $m$ and $\Ebar$ this has for
consequence that \eqref{shorttimeas} is indeed valid at least on this
time interval thus proving Theorem \ref{bornes}.

\section{Preservation of $\|\mu_N\|_{\infty,\eta}$}
From the form of the estimate on $m$ in Lemma \ref{boundm}, it is
clear that with this estimate we will never get a result for a
long time. Indeed, even assuming that we have bounded before $K$
and $R$, we would have the equivalent of $\dot m\leq m\times
\dEbar\leq C\,m\times m^{2+2\alpha'}$.

On the other hand this is somewhat strange since, in the limit,
the $L^\infty$ norm is conserved. And this preservation is very
useful in the proof of the existence and uniqueness of the
solution of the Vlasov equation, see for instance \cite{LP2}. But,
how to obtain the analog of this in the discrete case? At this
time, we just have a bound on $\|\mu_N\|_{\infty,\ep}$ on a small
time, and the bound is too huge to allows us to prove convergence
results for long time. Of course, this norm is not preserved at
all because we are looking at the scheme at the scale of the the
discretization. And in our calculation we do not use the fact that
the flow is divergence free, a property that is the key for the
preservation of the $L^\infty$ norm.

So what else can we do? One of the solution is to look at a scale
$\eta > \ep$, with $\ep/\eta$ going to zero as $\ep$ goes to zero.
At this scale, we have many more particles in a cell and we will
be able to obtain the asymptotical preservation of this norm. This
will be very useful because it will allow us to sharpen our
estimate on $E$ and $\delta E$. And with this we will obtain long
time convergence results.

\begin{figure}
\epsfxsize=13cm
\epsfysize=10cm
\epsfbox{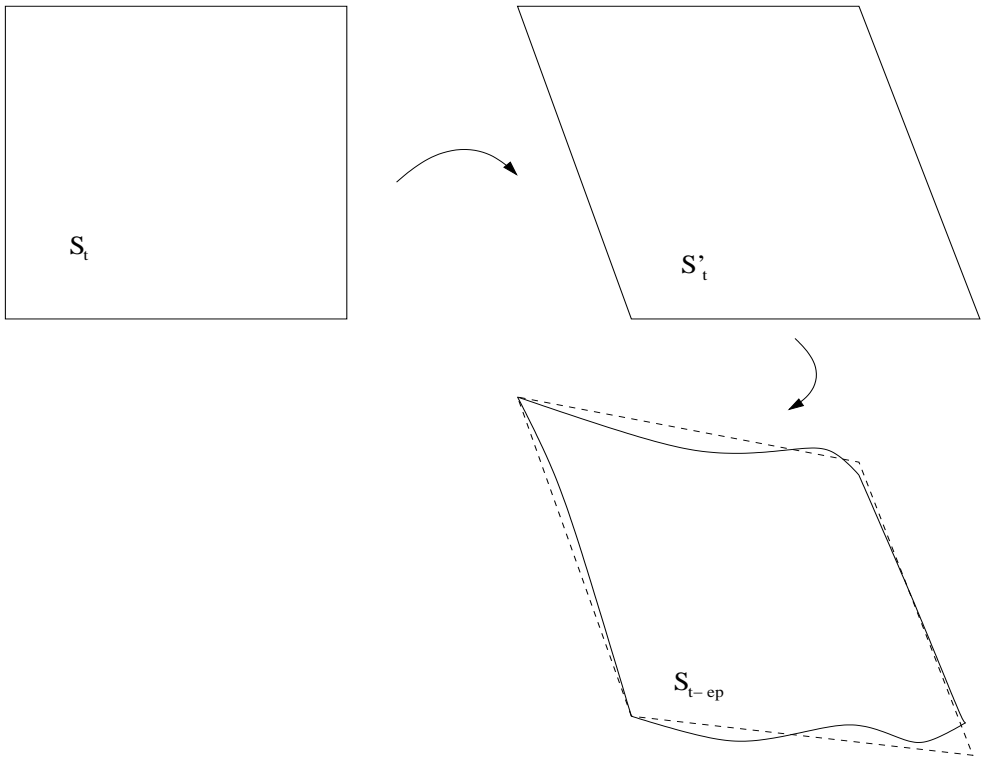}
\caption{Evolution of $S_t$.}\label{parallelo}
\end{figure}

Now, we will try to give roughly the idea of the proof in
dimension $1$ before beginning the genuine calculations. We choose
a time $t$ and a particle $i$, and look at the square of size
$\eta$ in the phase space centered on the particle $i$. We called
it $S_t= \{(x,v)| |x-X_i(t)|<\eta, |v-V(t)|<\eta\}$. We want an
estimation of the number of particles in this square at the time
$t$. For this, we first wanted to know where were these particles
at the time $t-\ep$. During the interval of time $[t-\ep,t]$, the
particle $j$ has moved of
\begin{eqnarray}
X_j(t)-X_j(t-\ep) & \thickapprox & \ep V_j(t) \\
V_j(t)-V_j(t-\ep) & \thickapprox & \ep E(t, X_j(t))
\end{eqnarray}
To make it correct, we will have to replace the second left hand
side by an average on an interval of size $\ep$ in time, as usual,
but we keep this for explanation. If we take the values in the
right hand side for exact, the particles we are interested in were
in the following domain $S_{t- \ep}$.
$$
S_{t- \ep} = \{(x,v)| |x-X_i(t)-\ep v| \leq \eta \quad \text{and}
\quad |v-V_i(t) -\ep E(t,x -\ep v )| \leq \eta\}
$$
There is two steps to obtain $S_{t-\ep}$ from $S_t$. First,
translate for a fixed $\bar{v}$, the set $S_t \cap \{v=\bar{v}\}$
of $\bar{v}$, for all $v$. We call the resulting set $S_t'$. Then,
for a fixed $\bar{x}$, translate $S_t' \cap \{x=\bar{x}\}$ of
$-\Ebar(x)$, for all $x$. With these steps, we see that $S_t$ and
$S_{t-\ep}$ have the same volume (see figure \ref{parallelo}).
But, if we keep the set $S_t'$, and do this another time, we will
obtain a set with a more strange shape and so one. So, we will
approximate $S_{t-\ep}$ by a parallelogram and iterate this step
from the parallelogram. We will choose the parallelogram
\begin{equation}\begin{split}
\bar{S}_{t-\ep} = \{(x,v)\;|\ &|x - X_i(t)-\ep v|<\eta\\
\text{and} \quad &|v-V_i(t)- \ep(x-X_i(t)) \nabla E_\ep(t,X_i(t))
-\ep \bar{E}(X_i(t))| \leq \eta \},
\end{split}
\end{equation}
where $E_\ep$ is an approximation at $\ep$ of the field, defined
by:
$$ E_\ep(x) = \frac{1}{N} \sum_{i \neq j}
\frac{x-X_i(t)}{(|x-X_i(t)|+\ep)^{(\alpha+1)}}.$$
We use this approximation to obtain a usable value of $\nabla E$.
This set has almost the same volume that $S_t$. And if we begin
with a parallelogram, we still obtain a parallelogram. What we
have to check is that this approximation make sense, that means
that we forgot or added a non relevant set of particles. For this,
we need that the parallelogram do not become too stretched, in
other words that the two sides do not become parallel on the
figure in dimension one. Remember that we are not exactly
interested by the volume of the set $S_0$ (the set that we obtain
iterating the process till we reached the initial time), but by
the number of vertices of the initial networking $\ep \Z^{2d}$
inside the parallelogram. The volume of the parallelogram is a
good approximation of this if its width is big with respect to
$\ep$.

We avoid this if the slope of the sides of the parallelogram are
distinct. This is ensured if $t$ and $\ep \sum_{k \leq t/\ep}
|\nabla E_\ep(kt/\ep,X_i(kt/\ep))|$ are smaller than $1/2$,
because this two quantities are respectively the tangent and
the cotangent of the angle between the sides of the parallelogram
and the $x$-axis. Remark that the sum is bounded by $t \dEbar(t)$.

This limits us in time, so we will obtain the conservation of the
$\|\mu_N\|_{\infty,\eta}$ only on a short time. But this is not a
problem. Call this time $T'$. After this, we choose a new scale
$\eta'$ sufficiently large, by instance $\sqrt{\eta}$. Since we
have proved the asymptotic preservation of
$\|\mu_N\|_{\infty,\eta}$, we can do the same think again,
replacing $\ep$ by $\eta$ and $\eta$ by $\eta'$. And obtain the
preservation of $\|\mu_N\|_{\infty,\eta'}$ on a new small interval
of time and so on.

Now, this is time to do the calculation for the proof of what we
claimed before. We just need to remark that we draw our pictures in
dimension one. In this case, we deal with true parallelogram. But we
will only prove result for $d \geq 2$ and in this case we will not use
true parallelograms. We will deal with sets defined by
$$ S=\left\{(x,v)| \left\|M \left(\begin{array}{c} x \\v \end{array}\right)
\right\| \leq \eta \right\} $$
where $M$ is a matrix of $\mathcal{M}_{2d}(\R)$, and the norm
$\|\cdot\|$ is defined by
$$
\|(x,v)\| = \max(|x|,|v|)
$$
For convenience, we will often decompose the matrix $M$ in four
blocks like below
$$
M= \left(\begin{array}{cc} A & B \\ C & D \end{array}\right)
$$

\begin{prop} \label{defpar}
Choose $(X_0,V_0) \in \R^{2dn}$, and $A$,$B$,$C$,$D$ four square
matrices in $\mathcal{M}_n(\R)$ satisfying the following conditions
that will be called the norm conditions in the rest of the paper
\begin{equation}
\label{normcond}
\max(\|A- Id\|,\|D-Id\|) \leq 1/2 \quad
\text{and} \quad \max(\|C\|,\|D\|) \leq 1/2,
\end{equation}
with the  dual norm on the matrix. We define $M$ as above and
denote by $N_t$ the number of particles in
$$
S_t=\left\lbrace (x,v)| \left\| M \cdot \left(\begin{array}{c} x-X_0 \\v-V_0
\end{array} \right) \right\|\leq \eta \right\rbrace
$$
We assume as in all the preceding results, that
$$
m(t_0) \leq \frac{1}{12\ep \,K(t)\dEbar(t)}
$$
Then, there exist a constant
$C=C(R,K,\Ebar,\|\mu\|_{\infty,\ep})$, there exists a position
$X'$ a speed $V'$, four matrices $A'$,$B'$,$C'$,$D'$ and then a
matrix $M'$ defined with this four blocks such that if we denote
$N_{t-\ep}$ the number of particles in
$$
S_{t-\ep}=\left\lbrace (x,v)| \left\| M' \cdot
\left(\begin{array}{c} x-X'_0 \\v-V'_0 \end{array} \right)
\right\|\leq \eta +C\ep(\eta^{\beta} +\ep) \right\rbrace
$$
we have the following inequality
$$
N_t \leq N_{t-\ep}.
$$
Moreover, we have
$$
\max(\|D-D'\|,\|D-D'\|,\|C-C'\|,\|D-D'\|) \leq C \ep.
$$
\end{prop}

\begin{proof} We divide it in two steps.

{\em Step 1: Estimate on $V_j(t-\ep) - v- \ep E(X_j(t))$.} Choose
 $j \in \{1,\dots,N\}$ such that $X_j(t) \in S_t$. We will first
work on the velocities and then integrate our estimation to obtain
what we need on the positions. We have
\begin{eqnarray*}
V_j(t-\ep) -V_j(t)& = & \ep \int_0^1 E(X_j(t-s\ep))\,ds \\
 & = & \ep \int_0^1 E(X_j(t-s\ep)) - E_\ep(X_j(t-s\ep))\,ds \\
 & & \qquad + \ep \int_0^1 E_\ep(X_j(t-s\ep))\,ds \\
 & = & I \quad +\quad II.
\end{eqnarray*}
We need to bound the first term $I$. The approximation
error is
\begin{multline*}
E(X_j(t-s\ep)) - E_\ep(X_j(t-s\ep))= \frac{1}{N} \sum_{k \neq j}
\Bigl(\frac{1}{|X_j(t')-X_k(t')|^{1+\alpha}} \\
-\frac{1}{(|X_j(t')-
X_k(t')|+\ep)^{1+\alpha}}\Bigr) (X_j(t')-X_k(t')).
\end{multline*}

To compute this term, we use again the same decomposition of the
phase space as in lemma \ref{moytpsdE}. If $k$ is such that $|X_k
- X_j| \geq \ep$, we can bound a difference in the sum by $
2^{1+\alpha}\ep / (|X_k-X_j|)^{1+\alpha} $ and then compute it
dividing the phase space in diadic subset. We obtain that
$$
\sum_{|X_k-X_j|\geq \ep} 2^{1+\alpha}\ep / (|X_k-X_j|)^\alpha
\leq C \ep.
$$
The constant $C$ that we obtain here depends of
$R$, $K$, $\Ebar$, $\|\mu_N\|_{\infty,\ep}$ exactly as in lemma
\ref{moytpsdE}. For the others terms, those were $|X_k - X_j| \leq
\ep$, we bound the difference by the sum of the norm of the two
term and do the same estimation as in (\ref{moytpsdE}). It is
possible because we do a mean over a interval of time of size
$\ep$. We obtain
$$
|I| \leq C\ep^2 + C'\ep^{d+1-\alpha}.
$$
again with the same dependence for $C'$. In the following, we will
assume to simplify the presentation that $d\geq 2$.  In that case,
we may replace it by $I \leq C\ep^2$ and this bound will be
enough. But all is still true if $d=1$. In that case, we replace
the last estimate by $I \leq C\ep^{2-alpha}$. Since $é-\alpha >
1$, all our estimates will remain valid.

Now, we approximate the second term $II$ by the same
expression where we replace $E_\ep$ by its first order
linearization near $X_0$. We defined $E_\ep^{lin}$ by
$$
E_\ep^{lin}(t,x) = E_\ep(t,X_0) - \nabla E_\ep(t,X_0) \cdot(x-X_0).
$$
We have to estimate the difference between $E_\ep$ and
$E_\ep^{lin}$. For this, we need to bound a sum of terms of the
following type
\begin{multline}
\frac{x-X_k}{(|x-X_k|+\ep)^{1+\alpha}} -
\frac{X_0-X_k}{(|X_0-X_k|+\ep)^{1+\alpha}} \\ - \nabla \left(
\frac{x-X_k}{(|x-X_k|+\ep)^{1+\alpha}}\right)(t,X_0) \cdot (x-X_0).
\end{multline}
For this, we find a path $I(s)$ between $x$ and $X_0$ (in others
words $I(0)=x$ and $I(1)=X_0$) of length smaller than $4|x-X_0|$
and such that $|I(t)- X_k| \geq \min(|x-X_k|,|X_0-X_k|)$. The
term to estimate may be rewritten
\[ \begin{split}
\int_0^1 \Big( \nabla
\Big( &\frac{x-X_k}{(|x-X_k|+\ep)^{1+\alpha}} \Big)(t,I(s)) \\
& -\nabla \left(
\frac{x-X_k}{(|x-X_k|+\ep)^{1+\alpha}}\right)(t,X_0) \Big) \cdot
I'(s) \,ds,
\end{split}\]
and each difference (without the multiplication by $I'(s)$) may be
bounded by
$$
C \frac{|I(t)-X_0|}{\min(|x-X_k|,|X_0-X_k|)^{2+\alpha}},
$$
where $C$ is just a numerical constant and also by
$$
\frac{C}{\min(|x-X_k|,|X_0-X_k|)^{1+\alpha}},
$$
where the constant $C$ is again independent of the problem. Doing a classical
interpolation between this two terms, we may bound it for every
$\gamma \in (0,1)$ by
$$
C
\frac{|I(t)-X_0|^{1+\gamma}}{\min(|x-X_k|,|X_0-X_k|)^{1+\alpha+\gamma}}.
$$
Thus, we may bound the approximation error by
$$
\frac{C}{N} \sum_{k \neq j}
\frac{|x-X_k|^{1+\gamma}}{\min(|x-X_k|,|X_0-X_k|)^{1+\alpha+\gamma}}.
$$
Now, we do the same computation as in the estimation of $\Delta
\bar{E}$. We obtain, as in lemma (\ref{moytpsdE}), that for every
$\gamma=\beta-1$, there exists a constant $C$ depending always on
$R$, $K$, $\Ebar$, $\|\mu_N\|_{\infty,\ep}$ such that
$$
|\frac{1}{\ep} \int_{t-\ep}^t E_\ep(s,X_j(t-s\ep))-
E_\ep^{lin}(s,X_j(t-s\ep)) \,ds| \leq C (|X_jj(t)-X_0|^{\beta} +
\ep).
$$
We approximate the term $II$ by
$$
|II - \ep \int_0^1 E_\ep^{lin}(s,X_j(t-s\ep))|\,ds \leq  C \ep
(|X_j(t-s\ep)-X_0|^{\beta}+ \ep).
$$
Moreover, we know that the approximated force field is Lipschitz
with a constant depending on $R$, $K$, $\Ebar$ and
$\|\mu_N\|_{\infty,\ep}$ on our small interval of time. This can
be seen just by doing the same calculation that gives the bound on
$\Delta \bar{E}$. Indeed, it is the same because
regularization at order $\ep$ or estimation on $\Delta E$ at order
$\ep$ are similar. Moreover, we may replace the value of
$E_\ep^{lin}$ taken at $X_j(t-s\ep)$ by those taken at $X_j(t)$.
Because of the bound on the speed, this will only introduce an
error of order $\ep$. Actually, if we denote $\tilde E = 1/\ep
\int_{t-\ep}^t E_\ep(s,X_0) \,ds$ and $\nabla \tilde{E } = 1/\ep
\int_{t-\ep}^t \nabla E_\ep(s,X_0) \,ds$, we have the following
inequalities
\begin{eqnarray}
|V_j(t-\ep) - V_j(t) - \ep (\tilde{E} +\nabla \tilde{E} \cdot
X_j(t-\ep)) | & \leq & C\ep(|X_j(t)-X_0| +\ep),\qquad \\
|X_j(t-\ep) - X_j(t) - \ep V_j(t-\ep)| & \leq & C \ep^2,
\end{eqnarray}
where the constant $C$ has the same dependance as before.

\medskip

{\em Step 2: The new parallelogram.} Now, we obtain an
approximated parallelogram that contains at time $t-\ep$ almost
all the particles that are in $S_t$ at time $t$. We look first at
the two first blocks of the matrix and  define the following
linear mappings
$$ A'= A + \ep  B \cdot \nabla {E_\ep}, $$
$$ B' = B + \ep A.$$
The center of our new parallelogram is given by
$$
X'_0 = X_0 - \ep V_0 \quad , \quad V'_0= V_0 - \int_{t-\ep}^t
E_\ep(s,X_0)\,ds.
$$
The previous calculation tells us that
$$
|A'\cdot (X_j(t-\ep) - X_0') - B' \cdot (V_j(t-\ep)- V_0') | \leq \eta
+ C\ep( |X_j(t-\ep)-X_0'|^{\beta} + \ep).
$$
Of course we obtain the same for the second line, and if we define also
$ C'= C + \ep  D \cdot \nabla {E_\ep} $, $ D' = D +\ep C$ and the
associated $M'$,
we obtain
$$
\left\| M' \left(\begin{array}{c} X_j(t-\ep)-X'_0 \\V_j(t-\ep)-V'_0
\end{array} \right) \right\|\leq C\ep(\|X_j(t) - X_0,V_j(t)-V_0\|^{\beta}
+ \ep).
$$
To conclude, we just have to prove that $|X_j(t) - X_0|$ is of
order $\eta$. But, this is true because our initial parallelogram
is not too stretched thanks to the conditions on the norm of
$A$, $B$, $C$, $D$. The following lemma is more precise.
\begin{lemma}
Assume that the four matrices $A$, $B$, $C$, $D$ satisfy
$\max(\|A-Id\|,\|D-Id\|,\|C\|,\|D\|) \leq 1/2 $. We use here the
dual norm on matrices, $\|A\|= \sup_{|X| \leq 1} |A \cdot X| /
|X|$. We denote as usual by $M$ the matrix composed of the four
blocks $A,\ B,\ C,\ D$, and
$$
S = \{(x,v)| \|M \cdot (x,v)^T \|
\leq \eta\}.
$$
Then
$$
S \subset B(0,2\eta).
$$
\end{lemma}
\begin{proof}[Proof of the lemma]
Choose an $x \in S$. With the definition of our norm $\|\cdot\|$
it means that
\begin{eqnarray*}
|A \,x + B\, v| & \leq & \eta, \\
|D \,v + C\, x| & \leq & \eta. \\
\end{eqnarray*}
This implies
\begin{eqnarray*}
|A \,x| - |B\, v| & \leq & \eta, \\
|D \,v| - |C\, x| & \leq & \eta. \\
\end{eqnarray*}
With the assumption on the matrix, we obtained
\begin{eqnarray*}
|x| - \frac{1}{2}(|x|+|v|) & \leq & \eta, \\
|v| - \frac{1}{2}(|x|+|v|) & \leq & \eta. \\
\end{eqnarray*}
We add this two lines and obtain
$$ |x| + |v| \leq 2 \eta. $$
\end{proof}
There only remains to prove the estimate on the determinant of
$M'$. But,
$$ M'-= \left(\begin{array}{cc} A + \ep B \, \nabla \tilde{E} & B +
\ep A\\ C +\ep D
\cdot \nabla \tilde{E}& D +\ep C\end{array}\right). $$
Then,
\begin{eqnarray}
\det(M') & = & \det \left(\begin{array}{cc} A + \ep B \, \nabla \tilde{E}
 & B + \ep A\\ C +\ep D  \nabla \tilde{E}& D +\ep C\end{array}\right) \\
 & =& \det \left(\begin{array}{cc} A - \ep^2 A \nabla \tilde{E}
& B + \ep A\\ C -\ep^2 C \nabla \tilde{E}& D +\ep C\end{array}\right) \\
\end{eqnarray}
To obtain the second line, we substract the second column
multiplied by $\nabla \tilde{E}$ to the first. Then, we see that
the new determinant is $\ep^2$ close from the preceding. This
gives us the expected bound on the determinant of $M'$.
\end{proof}
With the help of the proposition, we prove the almost
preservation of the $L^\infty$ norm at the scale $\eta$, which is
stated in the theorem below.
\begin{theorem} \label{preservlinfty}
There exist a time $T'$, a constant $C$ both depending on
$R(T')$, $K(T')$, $\Ebar(T')$, $sup_{t \leq T'}
\|\mu_N(t)\|_{\infty,\ep})$ such that if $t \leq T'$ and
\[
m(T')\leq \frac{1}{12\,\ep\,K(T')\,\dEbar(T')},
\]
the
following inequality holds:
$$
\|\mu_N(t)\|_{\infty, \eta} \leq \|\mu\|_{\infty,\ep} + C \,
(\eta^\beta + \ep/\eta).
$$
\end{theorem}
\begin{proof}
Let us choose a box $S_t = \{(x,v)\,|\;\|x-X_0,v-V_0\| \leq \eta\}$. We
can find a parallelogram containing at time $t-\ep$ the particles
that are in $S_t$ at time $t$. And we can iterate this process
till the parallelogram do not become too much stretched. In others
words, till the conditions \ref{normcond} are satisfied by the
matrix $M$. As at each step, the matrix move from at worst
$C'\ep$, we may iterate the process on an interval of time of
length $1/2C'$. So we obtain a parallelogram $S_0$ not too
stretched, if $t \leq 1/2C'$. We made $([t/\ep]+1)$ steps to reach
the time $0$. The last is necessary of length less than $\ep$, but
as in the proof of the theorem one, this raises no difficulty. So,
we know that $S_0$ is not too stretched, but we have to check that
this parallelogram $S_0$ is not too big.
\medskip

{\em Step 1: The size of $S_0$} We define real functions $g_\ep$
by $g_\ep(\eta)= \eta + C \ep (\ep + \eta^{\beta})$. We introduce
it because if $S_{t'}$ is defined by
\[
S_{t'}=\left\lbrace (x,v)| \left\| M \cdot \left(\begin{array}{c}
x-X_0 \\v-V_0
\end{array} \right) \right\|\leq \eta \right\rbrace.
\]
Then $S_{t'-\ep}$ is defined by a similar condition with new $M$,
$X_0$,$V_0$ where the right hand side of the inequality is
replaced by $g_\ep(\eta)$. And, so $S_0$ is defined by a similar
conditions, where the right hand side $\eta$ is replaced by
\[
\eta'=g_\ep \circ \dots \circ g_\ep(\eta),
\]
where the dots means $[t/\ep] + 1$ times. We need to control this
quantity. We define $r_n= g_\ep^n(\eta)$ where the exponent means
composed $n$ times. This $r_n$ is the quantity in the right
hand side of the definition of $S_{t-n\ep}$. From the formula
$r_{n+1}=g(r_n)$ we expect that $r_n \approx \eta + C n \ep (\ep +
\eta^{\beta})$. To prove this rigourously, we introduced
$\alpha_n$ defined by $r_n= \eta + Cn\ep(\ep + \eta^{\beta}) +
\alpha_n$. Provided $\alpha_n<\eta$, it
satisfies the following relation:
\[
\alpha_{n+1} \leq (1+C'\ep \eta^{\beta-1}) \alpha_n + C\,C'\ep^2
n\eta^{\beta-1}(\ep + \eta^{\beta}),
\]
where $C'$ is a numerical constant. We can bound $n\ep$ by $t$
since we will only iterate the process till $n= [t/\ep] +1$. The
previous inequality becomes
\[
\alpha_{n+1} \leq (1+C'\ep \eta^{\beta-1}) \alpha_n + C\,C'\ep\,
\eta^{\beta-1}(\ep + \eta^{\beta}).
\]
Since $\alpha_0=0$, we obtain that
\[
\alpha_n \leq ((1+C'\ep\eta^{\beta-1})^n -1)C(\ep + \eta^{\beta}).
\]
This gives us
\begin{eqnarray}
\eta' & \leq & \eta + (1 + e^{C'\eta^{\beta-1}} - 1) C (\ep +
\eta^{\beta}) \\
& \leq & \eta + C (\ep + \eta^{\beta}),
\end{eqnarray}
with a new constant for the last line.

\medskip

{\em Step 2: Covering of $S_0$ by balls of radius $\ep/2$}.
\begin{lemma} \label{cover}
let $S$ be a parallelogram defined as above by
$$
S=\left\lbrace (x,v)| \left\| M \cdot \left(\begin{array}{c} x-X_0 \\v-V_0
\end{array} \right) \right\|\leq \eta \right\rbrace
$$
with $M$ composed of the four blocks $A$, $B$, $C$, $D$ satisfying
the assumption \ref{normcond} and also $\det(M) \leq 1 + C\ep$.
Then, there exists a constant $C'$ depending only on $C$ such that
$S$ can be covered by $ [ \ep^{-2d} (Vol(S) + C \ep \eta^{2d-1})]$
balls of size $\ep$. In others words, there exists a finite set
$P$ of cardinal less than $\ep^{-2d}(Vol(S) + C'\ep \eta^{2d-1})$
such that $S \subset \bigcup_{p \in P} B(p,\ep)$
\end{lemma}

\begin{proof} [Proof of the lemma]
We define
$$S^+_{2\ep} = \left\lbrace (x,v)| \left\| M \cdot
\left(\begin{array}{c} x-X_0
\\v-V_0 \end{array} \right) \right\|\leq \eta + 2\ep \right\rbrace \, ,$$
$$S^+_{4\ep} = \left\lbrace (x,v)| \left\| M \cdot
\left(\begin{array}{c} x-X_0
\\v-V_0 \end{array} \right) \right\|\leq \eta + 4\ep \right\rbrace, $$
and $P = \ep \Z \cap S^+_{2\ep}$. We look at the set $P_\ep$
consisting of the union of all the balls (for the norm
$\max(|x|,|v|)$) of size $\ep$ centered at points of $P$, that is
$P_\ep = P + B(0,\ep)$. We will show that this set is included in
$S^+_{4\ep}$. For this, we choose $(x,v) \in P_\ep$. We associate
to this point the couple $(m,n)$ such that $\|(x-\ep m,v-\ep n) \|
\leq \ep$. Then,
\begin{eqnarray*}
\left\| M \cdot \left(\begin{array}{c} x \\v
\end{array} \right) \right\|  & \leq &
\left\| M \cdot \left(\begin{array}{c} x- \ep m\\v-\ep n
\end{array} \right) \right\| + \left\| M \cdot \left(\begin{array}{c} \ep m
\\\ep n \end{array} \right) \right\| \\
& \leq & \|M\| \ep + \eta + 2\ep \\
& \leq & \eta + 4\ep.
\end{eqnarray*}
In the last line, we use $\|M\| \leq 2$. This inequalities comes
from the condition \ref{normcond}. Therefore we have the inclusion
$P_\ep \subset S^+_{4\ep}$.

Moreover, if we choose a point $(x,v) \in S$, we can find a point
$(\ep m, \ep n)$of $\ep \Z^{2d}$ such that $\|(x-\ep m,v-\ep n) \|
\leq \ep$. As above, we have
$$ \left\| M \cdot \left(\begin{array}{c} \ep m \\ \ep n
\end{array} \right) \right\| \leq \eta + 2\ep.$$
Thus, $\ep(m,n) \in P$. That proves that $S \subset P_\ep$. So, we
have the inclusions
$$ S \subset P_\ep \subset S_{4\ep}^+.$$
The first is the recovering we want. The second gives us an
estimate on the cardinal of $P$. Comparing the volume of $P_\ep$
and $S^+_{4\ep}$ we obtain
$$ (2\ep)^{2d} |P| \leq det(M)(\eta + 4\ep)^{2d}.$$
If $M$ satisfies $\det(M) \leq 1 + C\ep$, that gives us
$$ |P| \leq \ep^{-2d}\eta^{2d-1}(\eta + C'\ep).$$
\end{proof}

\medskip

{\em Step 3: Conclusion of the proof.} We choose $t$ as in the
previous step and a box $S_t$ in the phase space of size
$\eta$. We define the parallelogram $S_0$ as above. It
exits and is not too stretched because we choose $t \leq T'$. From
Step 1, we know that
\begin{eqnarray}
Vol(S_0) & \leq & (1+C\ep)(\eta + C(\ep + \eta^\beta))^{2d} \\
& \leq & \eta ^{2d}(1+ C(\eta^\beta +\ep)).
\end{eqnarray}
And since it is not too stretched, we may used lemma \ref{cover}
and cover it by less than $\ep^{-2d}(Vol(S_0) + C'\ep
\eta^{2d-1})$ balls of radius $\ep$. Actually, that means that
$$
N_t \leq N_0 \leq \|\mu_N(0)\|_{\ep,\infty} \eta^{2d}(1+
C(\eta^\beta +\ep),
$$
and dividing by $\eta^{2d}$ we obtain
$$
\|\mu_N(t)\|_{\eta,\infty} \leq \|\mu_N(0)\|_{\ep,\infty} (1+
C(\eta^\beta +\ep)).
$$
\end{proof}

\subsection{New estimates on $\Ebar$ and $\dEbar$}
%
 The almost preservation of the $\|\mu_N\|_{\infty,\eta}$ norms will
enable us to prove
a new estimate on $\Ebar$. We can obtain it by doing the same
separation of the position space in dyadic cells, but we begin
with cells $\tilde C_k$ satisfying
\[
C_k=\biggl\{i\,\Bigl|\; 3\,\eta\,K(t_0)\,
2^{k-1}<|X_i(t_1)-X_1(t_1)|\leq 3\,\eta\,K(t_0)\, 2^{k}\biggr\}.
\]
That gives a first term in $\|\mu_N\|_{\infty,\eta} K^{\alpha'}
R^{\alpha'-\alpha}$. Next we decompose $\tilde C_0$ in a union of
$C_k$ (the ``true'' cells with size $\ep$); Since the index $k$
goes only up to $\log(\eta/\ep)/\log 2$, the corresponding term is
less than $\|\mu_N\|_{\infty,\ep}\,
K^{2\alpha'-\alpha}\,\eta^{\alpha'-\alpha}$. The remainder term
coming for $C_0$ is dealt just as in Lemma \ref{moytps} and it
yields the two same terms.

That gives us the following lemma
\begin{lemma}  For any $\alpha'$ with $\alpha<\alpha'<3$, assume that
\[
m(t_0)\leq \frac{1}{12\,\ep\,K(t_0)\,\dEbar(t_0)},
\]
then the following inequality holds
\[ \begin{split}
\Ebar(t_0) \leq C\, ( &\|\mu_N\|_{\infty,\eta}^{\alpha'/d}\,
K^{\alpha'}\,R^{\alpha'-\alpha} +
\|\mu_N\|_{\infty,\ep}^{\alpha'/d}\,K^{2\alpha'-\alpha}\,
\eta^{\alpha'-\alpha}\\
&\ +\ep^{d-\alpha}\, \|\mu_N\|_{\infty,\ep}\, K^{2d-\alpha}
 + \ep^{2d-3\alpha}\,
\|\mu_N\|_{\infty,\ep}\, K^{d-\alpha} \Ebar^d),
\end{split}
\]
where we use the values of $\|\mu_N\|_{\infty,\ep}$, $R$, $K$, $m$
and $\Ebar$ at the time $t_0$. \label{newestEbar}
 \end{lemma}
The only non-negligable term in this estimate is sub-linear if
$\alpha'$ is chosen sufficiently close to $\alpha$.

Of course we can perform the same changes in the proof of $\dEbar$
to get
\begin{lemma}  For any $\alpha'$ with $\alpha<\alpha'<3$, assume that
\[
m(t_0)\leq \frac{1}{12\,\ep\,K(t_0)\,\dEbar(t_0)},
\]
then the following inequality holds
\[ \begin{split}
\dEbar(t_0) \leq C\, ( &\|\mu_N\|_{\infty,\eta}^{(1+\alpha')/d}\,
K^{1+\alpha'}\,R^{\alpha'-\alpha} +
\|\mu_N\|_{\infty,\ep}^{(1+\alpha')/d}\,K^{1+2\alpha'-\alpha}\,
\eta^{\alpha'-\alpha}\\
&\ +\ep^{d-\alpha-\beta}\, \|\mu_N\|_{\infty,\ep}\, K^{2d-\alpha}
 + \ep^{2d-3\alpha-\beta}\,
\|\mu_N\|_{\infty,\ep}\, K^{d-\alpha} \Ebar^d),
\end{split}
\]
where we use the values of $\|\mu_N\|_{\infty,\ep}$, $R$, $K$, $m$
and $\Ebar$ at the time $t_0$. \label{newestdEbar}
 \end{lemma}

\subsection{Proof of Theorem \ref{longtimeexist}}
Let us fix any time $T>0$. The aim is to show that we have bounds
for $R$, $K$, $\Ebar$ and $m$, uniform in $N$ on $[0,\ T]$.

Next we choose $\eta_0=\ep^{1/2}$ for instance and
$\eta'=\ep^{1/4}$.

Since for any $N$ the quantities $R$, $K$, $\Ebar$ and $m$ are
continuous in time, we may define $T_N<T$ as the first time $t$
(if it exists) such that one of the following inequality at least
is not true for some integer $M$ to be chosen after
\begin{equation}\begin{split}
&T'=T(R(t), K(t),\Ebar(t), \sup_{s\leq t}\|\mu_N\|_{\infty,\ep})\geq
\frac{T}{M},\\
&m(t)\leq \frac{1}{12\,\ep\,K(t)\,\dEbar(t)},\quad
C(R(t), K(t),\Ebar(t), \sup_{s\leq t}\|\mu_N\|_{\infty,\ep})
\leq \ep^{-1/8M},\\
&\ep^{d-\alpha} (m(t))^{-2d}\,(K(t))^{2d-\alpha}\leq
\ep^\beta,\quad
\ep^{2d-3\alpha}(m(t))^{-2d}\,(\Ebar(t))^d\,(K(t))^{d-\alpha}\leq
\ep^\beta.
\end{split}\label{ineq}\end{equation}
The quantity $T'$ and $C$ are the time and constant defined in Theorem
\ref{preservlinfty}.
Therefore on $[0,\ T_N]$ all inequalities
\eqref{ineq} are true and we may apply both Theorem 3.1 and Lemma
\ref{newestEbar}.

We define $t_i=i\,T'$ and $\eta_i=\eta_0\times r^i$ with
$r=\ep^{-1/4M}$ so that $\eta_M=\eta'$. We are going to apply $M$
times Theorem \ref{preservlinfty}, once on every interval
$[t_{i-1},\ t_{i}]$ (of size less than $T'$) and with
$\eta=\eta_i$ and $\ep$ replaced by $\eta_{i-1}$. That gives
\[
\sup_{t\in[t_{i-1},t_i]}\; \|\mu_N(t)\|_{\infty,\eta_i}\leq
\|\mu_N\|_{\infty,\eta_{i-1}}+C(\Ebar(t_i),\dEbar(t_i))\,
(\eta^\gamma_i+\ep^{1/4M}),
\]
and consequently thanks to \eqref{ineq}
\begin{equation}
\sup_{t\leq T_N} \|\mu_N(t)\|_{\infty,\eta'}\leq
\|\mu_N\|_{\infty,\ep}+C(\Ebar(T_N),\dEbar(T_N))\,M\,\ep^{1/4M}\leq
2\,\|\mu_N^0\|_{\infty,\ep}, \label{bornelinfty}
\end{equation}
independently of $N$ (and $T_N$). Now we apply Lemma
\ref{newestEbar} at time $T_N$ and because of \eqref{ineq}, we
obtain
\begin{equation}\begin{split}
\Ebar(T_N)&\leq C\,\|\mu_N(T_N)\|_{\infty,\eta}^{\alpha'/d}\,
(K(T_N))^{\alpha'}\,(R(T_N))^{\alpha'-\alpha}\\
&\leq C\,(K(T_N))^{\alpha'}\,(R(T_N))^{\alpha'-\alpha},
\end{split}\label{intboundEbartn}\end{equation}
using \eqref{bornelinfty}. As $T_N>\ep$, Lemma \ref{boundK}
implies that
\[
K(T_N)\leq K(0)+C\,\int_0^{T_N} \Ebar(t)\,dt\leq
K(0)+C\,T_N\,\Ebar(T_N).
\]
From this inequality, we immediately deduce that
\[
R(T_N)\leq R(0)+T_N\,K(0)+C\,T_N^2\,\Ebar(T_N)\leq
C\,T+C\,T^2\,\Ebar(T_N).
\]
Inserting these last two inequalities in \eqref{intboundEbartn},
we find
\[
\Ebar(T_N)\leq C\,T+C\,T^2\,(\Ebar(T_N))^{2\alpha'-\alpha}.
\]
Since $2\alpha'-\alpha<1$, there exists a constant $C(T)$
depending only on $T$ and the initial distribution such that
\begin{equation}
\Ebar(T_N)\leq C(T),\quad K(T_N)\leq C(T),\quad R(T_N)\leq C(T).
\label{longtimeboundEKR}
\end{equation}
We are almost ready to conclude, we only need to apply once Lemma
\ref{newestdEbar} and by \eqref{ineq}, \eqref{bornelinfty} and
\eqref{longtimeboundEKR}
\begin{equation}
\dEbar(T_N)\leq C(T). \label{longtimebounddE}
\end{equation}
Inserting \eqref{longtimebounddE} in Lemma \ref{boundm}, we
eventually get
\begin{equation}
m(T_N)\leq C(T). \label{longtimeboundm}
\end{equation}
Together \eqref{longtimeboundEKR}, \eqref{longtimebounddE} and
\eqref{longtimeboundm} imply that all the inequalities of
\eqref{ineq} are true with a factor $1/2$ at time $T_N$, provided
$N$ and $M$ are large enough. Therefore \eqref{ineq} is still true on at
least a short time interval after $T_N$ and that means that
necessarily $T_N=T$. The consequence is that
\eqref{longtimeboundEKR}, \eqref{longtimebounddE} and
\eqref{longtimeboundm} are true on any time interval $[0,\ T]$
which is exactly Theorem \ref{longtimeexist}.

Finally note that we have implicitly used the short time result
when we said that $T_N>\ep$.

\section{Convergence of the density in the approximation}
The existence of the bound on $R$, $K$, $\Ebar$, $\dEbar$ and $
\|\mu_N\|_{\infty,\eta}$ implies the weak convergence of the
distribution $\mu_N$ to a weak solution of the Vlasov equation and
Theorem \ref{derivation} is only a consequence of Theorem \ref{longtimeexist}
and the following proposition
\begin{prop}
Let $\mu_N$ be the distributions associated with the solutions to
\eqref{ODE}.
We assume that the initial conditions $\mu_N^0$ converges weakly in
$M^1(\R^{2d})$ to some
$f_0 \in L^1\cap\Linf(\R^{2d})$. We choose a time $T > 0$. Assume
furthermore that
 there exists a
constant $C(T)$ independent of $N$ such that
$$
\sup_{\ep>0}(R(T),\;K(T),\;\Ebar(T),\;\dEbar(T),\;\|\mu_{\infty,\eta}\|)
< +\infty \; ,
$$
where $\eta$ depends on $\ep$ and $N$ and goes to zero when $\ep$ goes to
zero. Then, $\mu_N(t)$ converges weakly to $f(t)$, a solution to
the Vlasov equation with initial conditions $f^0$.
\end{prop}
\begin{proof}
We recall that the distribution of the particles satisfies the
Vlasov equation in the sense of distribution. Moreover, the
sequence $\mu_N$ is bounded in $\CCC([0,T],M^1(\R^{3d}))$. Up to an
extraction, we may assume that $\mu_N$ converges weakly to some $f
\in \Linf([0,T],M^1(\R^{2d}))$. Moreover, the fact that
$\|\mu_N\|_{\infty,\eta}$ is bounded implies that $f \in \Linf$.
To see this, we choose a regular test function $\phi$ with compact
support. We have
$$
\langle \mu_N, \phi \rangle = \frac{1}{N} \sum_{i = 1}^N
\phi(X_i(t),V_i(t)).
$$
Now, we define $\rho_\eta(x,v)= \chi_C(x/\eta,v/\eta)$ where
$\chi_C$ is the characteristic function of the set $C =
\{(x,v)|\|(x,v)\| \leq 1\}$ and we write
\[ \begin{split} \langle \mu_N, \Phi \rangle
& = \frac{1}{N} \sum_{i = 1}^N \Phi
\ast \rho_\eta (X_i(t),V_i(t)) \\
& +\frac{1}{N} \sum_{i = 1}^N (\Phi(X_i(t),V_i(t)) - \Phi\ast
\rho_\eta (X_i(t),V_i(t))).
\end{split} \]
The first term is $\int \phi\ast\rho_\eta(x,v)\,d\mu_N(x,v) = \int
\phi (\mu_N*\rho_\eta)\,dxdv$. So it is bounded by $\|\phi\|_1
\|\phi \mu_N*\rho_\eta\|_{\infty}$. But $\|\phi
\mu_N*\rho_\eta\|_{\infty}$ is exactly $\|\mu_N\|_{\infty,\eta}$.
The second term is easily bounded by $\eta \|\nabla
\Phi\|_{\infty}$. Putting all together, we obtain that
$$
\langle \mu_N, \Phi \rangle \leq \|\mu_N\|_{\infty,\eta}
\|\Phi\|_1 + \eta \|\nabla \Phi\|_{\infty}.
$$
At the limit,
$$
\langle f, \Phi \rangle \leq \liminf_{N \rightarrow \infty}
\|\mu_N\|_{\infty,\eta} \|\Phi\|_1,
$$
which means that $f \in \Linf$ and that $\|f\|_\infty \leq
\liminf_{N \rightarrow \infty} \|\mu_N\|_{\infty,\eta}$.

\medskip

The passage to the limit in the linear part of the equation does not
raise any difficulty. For the term in $F \cdot \nabla_v f$, we
need a strong convergence in the force. We denote by $F_\infty$
the force induced by $f$ and by $F_N$ the force induced by
$\mu_N$
\begin{eqnarray*}
F_\infty (x) & = & \int \frac{x-y}{|x-y|^{1+\alpha}} \,dydw, \\
F_N(x) & = & \frac{1}{N} \sum_{i=1}^N
\frac{x-X_i(t)}{|x-X_i(t)|^{1+\alpha}}.
\end{eqnarray*}
We have
\[ \begin{split}
\frac{1}{\ep}\int_{t_0}^{t_0+\ep} F_N(X_i(t)) -&
F_\infty(X_i(t)) \,dt =I_1+I_2+I_3\\
&\quad = \frac{1}{\ep}\int_{t}^{t+\ep}\int_{|y-X_i(s)| \geq r}
  \frac{y-X_i(s)}{|y-X_i(s)|^{\alpha+1}} \,d(\mu_N - f)(y)\,ds\\
&\qquad + \frac{1}{\ep} \int_t^{t+\ep} \int_{|y-X_i(t)| \leq r}
  \frac{y-X_i(s)}{|y-X_i(s)|^{\alpha+1}} \,d\mu_N(y) \,ds  \\
&\qquad - \frac{1}{\ep}\int_{t}^{t+\ep}\int_{|y-X_i(t)| \leq r}
  \frac{y-X_i(s)}{|y-X_i(s)|^{\alpha+1}} \,df(y) \,ds,
\end{split}
\]
for all $r>0$. The first term $I_1$ in the right hand side always
goes to zero because $\mu_N$ converges weakly to $f$. The second
term is dominated by $\|f\|_\infty \int_{B(0,R)} dy/|y|^\alpha$, a
quantity which is less than $C \|f\|_\infty r^{d-\alpha}$. The
last one is the field created by the close particles in the
discrete case. To estimate it we will use the lemma \ref{moytps},
replacing $R(T)$ by $r$. This gives
\[ \begin{split}
 I_3 \leq C \, (&\|\mu_N\|_{\infty,\ep}^{\alpha'/d}\,
K^{\alpha'}\,r^{\alpha'-\alpha} + \ep^{d-\alpha}\,
\|\mu_N\|_{\infty,\ep}\,
K^{2d-\alpha} \\
&\ + \ep^{2d-3\alpha}\, \|\mu_N\|_{\infty,\ep}\, K^{d-\alpha}
\Ebar^d\,K^d)\leq C\, r^{\alpha'-\alpha}.
\end{split} \]
And these bounds are independent of $N$ or $i$.

Then, letting $\ep$ going to $0$ and then $r$, we find that
\begin{equation}
\label{Fconv} \sup_{i, t} \frac{1}{\ep} \int_{t}^{t + \ep}
|F_N(X_i(s)) - F_\infty(X_i(s))| \,ds \rightarrow 0 \quad
\text{as} \quad \ep \rightarrow 0.
\end{equation}
With this strong convergence, we are able to prove the convergence
of the term $F_N \cdot \nabla_v \mu_N \,ds$ towards $F_\infty \cdot
\nabla_v f$ in the sense of distributions. We
choose a test smooth test function $\phi$ with compact support and
compute
\begin{multline}
J = \int_0^T \Big (\int_{x,v} F_\infty(t,x) \cdot \nabla_v
\phi(t,x,v) f(t,x,v) \, dxdv  \\
- \sum_{i=1}^N F_N(t,X_i(t),V_i(t)) \cdot \nabla_v
\phi(t,X_i(t),V_i(t)) \Big) \,dt
\end{multline}
We separate $J$ in $J_1 + J_2$, with
$$
J_1 = \int_0^T  \int_{x,v} F_\infty(t,x) \cdot \nabla_v
\phi(t,x,v) d(f-\mu_N)(,x,v) \,dt,
$$
and
$$
J_2 = \int_0^T \Big( \sum_{i=1}^N F_\infty(t,X_i(t),V_i(t))-
F_N(t,X_i(t),V_i(t)) \cdot \nabla_v \phi(t,X_i(t),V_i(t)) \Big)
\,dt
$$
Because of the continuity of $F_\infty$, $J_1$ vanishes as
$\ep$ goes to zero. To show that $J_2$ vanishes as well, we decompose
it in $M= [T/\ep]+1$ integrals on $M$ intervals of time with length $\ep$.
The last interval is of length less than $\ep$, but that does not any
difficulty and we do as if it were of length $\ep$. We obtain,
\begin{equation}\begin{split}
J_2 & = \sum_{k=1}^M \int_{k\ep}^{(k+1)\ep} \Bigl(\sum_{i=1}^N
\bigl(F_\infty(t,X_i(t),V_i(t))- F_N(t,X_i(t),V_i(t))\bigr)\\
&\qquad\quad \cdot \nabla_v
\phi(t,X_i(t),V_i(t)) \Bigr) \,dt \\
& \leq  C \sum_{k=1}^M \int_{k\ep}^{(k+1)\ep} \Big(\sum_{i=1}^N
\bigl|F_\infty(t,X_i(t),V_i(t))- F_N(t,X_i(t),V_i(t))\bigr|\Big) \,dt.
\end{split}\end{equation}
This sum may be bounded by
$$
CT \sup_{i, t} \frac{1}{\ep} \int_{t}^{t + \ep} |F_N(X_i(s)) -
F_\infty(X_i(s))| \,ds,
$$
a quantity which goes to zero according to (\ref{Fconv}). Thus,
$J$ goes to zero when $\ep$ goes to zero and the proof is done.
\end{proof}

\end{document}